\date{} 
\title{KPZ in one dimensional random geometry\\of multiplicative cascades}
\author{%
Itai Benjamini\thanks{Microsoft Research and the Weizmann Institute of Science}
\and
Oded Schramm\thanks{Microsoft Research}
}
\documentclass[12pt]{article}
\usepackage{amsmath}
\usepackage{amssymb}
\usepackage{amsthm}
\usepackage{amsfonts}
\usepackage{graphicx}

\newif\iffigures\figurestrue

\newif\ifhyper\IfFileExists{hyperref.sty}{\hypertrue}{\hyperfalse}
\hyperfalse 

\ifhyper\usepackage[hypertex]{hyperref}
\def\hitem#1#2{\item[\hypertarget{#1}{#2.}]\expandafter\gdef\csname LBL#1ITM\endcsname{#2}}%
\def\iref#1{\hyperlink{#1}{\csname LBL#1ITM\endcsname}}%
\else
\def\hitem#1#2{\item[{#2.}]\expandafter\gdef\csname LBL#1ITM\endcsname{#2}}%
\def\iref#1{{\csname LBL#1ITM\endcsname}}%
\fi

\newif\ifdraft
\drafttrue 
\long\def\comment#1{}
\long\def\old#1{}

\numberwithin{equation}{section}
\numberwithin{figure}{section}

\newtheorem{theorem}{Theorem}
\numberwithin{theorem}{section}

\newtheorem{lemma}[theorem]{Lemma}

\theoremstyle{remark}
\theoremstyle{remark}
\newcounter{mycount}

\let\qqed=\qed
\def\QED{\qqed\medskip}
\let\qed=\QED

\newcommand{\Q}{\mathbb{Q}}

\newcommand{\N}{\mathbb{N}}

\def\SLEkk#1/{$\mathrm{SLE}(#1)$}
\def\SLEr#1/{$\mathrm{SLE(\kappa;#1)}$}
\def\SLEkr#1;#2/{$\mathrm{SLE(#1;#2)}$}
\def\SLEk/{\SLEkk{\kappa}/}
\def\SLEtwo/{\SLEkk2/}
\def\SLE/{$\mathrm{SLE}$}
\def\SLEab/{\SLEkr 4; {a/\hco-1}, {b/\hco-1}/}

\def\Ito/{It\^o}
\def \eps {\epsilon}
\def \P {{\bf P}}
\def\md{\mid}
\def\Bb#1#2{{\def\md{\bigm| }#1\bigl[#2\bigr]}}
\def\BB#1#2{{\def\md{\Bigm| }#1\Bigl[#2\Bigr]}}
\def\Bs#1#2{{\def\md{\mid}#1[#2]}}
\def\Pb{\Bb\P}
\def\Eb{\Bb\E}

\def\EB{\BB\E}

\def\Es{\Bs\E}

\def \E {{\bf E}}

\def\ev#1{{\mathcal{#1}}}
\def\II{\mathcal I}

\def\Daa{{\Delta_0}}
\def\Dbb{\Delta}
\def\Da{{\zeta_0}}
\def\Db{\zeta}
\def \proof {{ \medbreak \noindent {\bf Proof.} }}
\def\proofof#1{{ \medbreak \noindent {\bf Proof of #1.} }}

\def\Mattila{MR1333890}
\def\KahanePeyrier{MR0431355}

\def\noopsort#1{}

\begin{document}
\maketitle

\begin{abstract}
We prove a formula relating the Hausdorff dimension of a subset of the unit interval
and the Hausdorff dimension of the same set with respect to a random path matric on the interval,
which is generated using a multiplicative cascade.
When the random variables generating the cascade are
exponentials of Gaussians,
the well known KPZ formula of  Knizhnik, Polyakov and Zamolodchikov from quantum gravity~\cite{MR947880} appears.
This note was inspired by the recent work of Duplantier and Sheffield~\cite{DSinprep} proving a somewhat different
version of the KPZ formula for Liouville gravity.
In contrast with the Liouville gravity setting,
the one dimensional multiplicative cascade framework facilitates the determination
of the Hausdorff dimension, rather than some expected box count dimension.
\end{abstract}


\section{Introduction}
There is growing interest in establishing a rigorous theory of two dimensional continuum quantum gravity.
Heuristically, quantum gravity is a metric chosen on the sphere uniformly among all possible metrics.
Although there are successful discrete mathematical quantum gravity models, we do not yet
have a satisfactory continuum definition.
A highlight of quantum gravity in the physics literature is the mysterious KPZ formula of
Knizhnik, Polyakov and Zamolodchikov~\cite{MR913671, MR947880}
relating the dimensions of fractals
in the random geometry to the corresponding dimension in Euclidean geometry.
More specifically, the KPZ formula is
\begin{equation}
\label{e.KPZ}
\Dbb-\Daa=\frac{\Dbb(1-\Dbb)}{k+2}\,,
\end{equation}
where $2-2\,\Dbb$ is the dimension of a set in quantum gravity and
$2-2\,\Daa$ is the dimension of the
corresponding set in Euclidean geometry\footnote{The parameter $k$ comes up since it is presumed that
there is essentially one free parameter in the construction of quantum gravity.
This parameter is intimately related to the central charge, and the various variants of
quantum gravity are believed to arise by weighting a uniform measure with the partition
function of a statistical physics model.}.

Recently Duplantier and Sheffield~\cite{DSinprep} were able to prove an expected box count dimension version of the KPZ formula
in Liouville gravity. In their setup, they avoid the very difficult issue of defining the random metric\footnote{We use the term \lq\lq metric\rq\rq\ to mean a \lq\lq distance function\rq\rq, not a Riemannian metric tensor.}, and instead define a random measure.

Multiplicative cascades is a well studied object and
defines naturally a random metric $\rho$ on $[0,1]$.
(The definition will be recalled in Section~\ref{s.bg}.)
The metric space $\bigl([0,1],\rho\bigr)$ is a path metric space, which in this
case just means that it is isometric with some interval $[0,\ell]$ with its Euclidean
metric $|x-y|$. The length $\ell=\rho(0,1)$ of $[0,1]$ in the metric $\rho$ is in general random.
In this note we  prove a formula
relating the Hausdorff dimensions of sets $K\subset [0,1]$ with respect to the random metric $\rho$ on the one hand and
with respect to the Euclidean metric on the other hand.
The KPZ relation appears precisely when the random variables defining the cascade are exponentials of
Gaussians.
Interestingly, this is hinted at in the introduction of~\cite{MR1835030}.
 
One goal of the note is to establish in this simple one-dimensional setup a Hausdorff dimension version of
the KPZ relation.  We also generalize the discussion to the non-Gaussian setting, mainly to gain
perspective on the underpinnings of the KPZ relation.

The proof in~\cite{DSinprep} is based on large deviations arguments, 
which is not the case for the present paper. While in a very general sense one could
say that the ideas are similar, it can reasonably be claimed that our proof is substantially
different and is rather elegant. As we have not yet seen a draft of~\cite{DSinprep},
this comparison is based on lectures by Scott Sheffield and on conversations with him.

\medskip
 
Early versions of multiplicative cascades were introduced by 
Kolmogorov already in 1941~\cite{MR1124922} and were developed by Yaglom~\cite{YaglomTurbulence} and by
Mandelbrot~\cite{MandelbrotCascades}.
Many fundamental properties of multiplicative cascades were first proved
in the remarkable paper~\cite{\KahanePeyrier} by Kahane and Peyri\`ere.
For further
background and references regarding the long history of multiplicative cascades see e.g.,~\cite{MR1835030}.

\medskip

In Section~\ref{s.bg}, we describe the basic setup of multiplicative cascades and
define the random metric $\rho$.
Our main result is stated and proved in Section~\ref{s.haus}.
An appendix follows in which we prove some
essentially known necessary background facts about
multiplicative cascades which we need to use. In some cases, the proofs
in the appendix are simpler than the proofs that we were able to find
in the literature, and in other cases, the results are slightly stronger.

\bigskip
\noindent
{\bf Acknowledgements}: We are obliged to Scott Sheffield for numerous discussions explaining his insights to us.
Yuval Peres and Ed Waymire have been very helpful in enlightening us
with regards to the published work on multiplicative cascades.

\section{Setup}\label{s.bg}

We now describe our setup.
Let $\II_n$ denote the set of dyadic subintervals of $[0,1]$ of length $2^{-n}$;
namely,
$$
\II_n:=\Bigl\{ \bigl[k\,2^{-n},(k+1)\,2^{-n}\bigr]: n\in\N,\,k\in \{0,1,\dots,2^{n}-1\}
\Bigr\}.
$$
Then each interval in $\II_n$ has precisely two subintervals in $\II_{n+1}$,
its left half and its right half.
Also set $\II =\bigcup_{n\in\N} \II_n$.

Let $W$ be some positive random variable with mean $1$,
and let $W_{I}$, $I\in\II$, be an independent collection
of random variables, each of which has the distribution of $W$.
We now define inductively a sequence of random measures on $[0,1]$.
Let $\mu_0$ denote Lebesgue measure on $[0,1]$,
and let $\mu_1:= W_{[0,1]}\,\mu_0$.
Let $\mu_2$ denote the measure which agrees with $W_{[0,1/2]}\,\mu_1$
on $[0,1/2]$ and agrees with $W_{[1/2,1]}\,\mu_1$ on $[1/2,1]$.
Inductively, define $\mu_{n+1}$ as the measure
that on every $I\in \II_{n}$ agrees with $W_I\,\mu_n$.
Alternatively, 
$$
\mu_{n}:= w_n\,\mu_0\,,
\qquad\text{ where }\qquad
w_n(x):=\prod_{j=0}^{n-1} W_{I_j(x)}\,,
$$
and $I_j(x)$ denotes the interval $I\in\II_j$ that contains $x$
(and if there is more than one, the one whose maximum is $x$, say).

We will need the following general result regarding multiplicative
cascades.

\begin{theorem}\label{t.bg}
The weak limit $\mu:=\lim_{n\to\infty}\mu_n$ exists almost surely.
Moreover, if $\Eb{W\log_2 W}<1$ then $\mu[0,1]>0$ a.s.\ and $\mu$
has no atoms a.s.
\end{theorem}

This theorem is entirely or almost entirely proved in~\cite{\KahanePeyrier}, but
we present in the appendix a different perhaps simpler proof.
Since $\mu[0,s]$ is a positive martingale for every $s\in[0,1]$, the first
claim is very easy to verify.
In~\cite{\KahanePeyrier} they show that
$\mu[0,1]>0$ a.s.\ if and only if $\Eb{W\log_2 W}<1$.
The claim about the non-existence of atoms would follow from the last remark
in~\cite{\KahanePeyrier}, but we were unable to verify the justification
of that remark (though it does follow under additional moment assumptions).

\smallskip

Henceforth, we will be assuming that
\begin{equation}
\label{e.WlogW}
\Eb{W\log_2W}<1\,.
\end{equation}
By the result of~\cite{\KahanePeyrier} mentioned above, this assumption is
necessary for the limit $\mu$ to be nonzero.

On $[0,1]$, define the random metric $\rho$ by
$$
\rho(x,y):=\mu[x,y]\,,\qquad \text{ for all }0\le x\le y\le 1\,.
$$
If $F(x):=\mu[0,x]$, then $\rho$ is just the pullback of the
Euclidean metric on $[0,F(1)]$ under $F$.
We also set
\begin{equation}
\label{e.ell}
\ell_n:=\mu_n[0,1]\,,\qquad\ell:=\mu[0,1]\,.
\end{equation}
Clearly, $\Eb{\ell_n}=1$ and $\Eb{\ell}\le 1$.
In fact, $\Eb{\ell}=1$ by~\cite{\KahanePeyrier}, but we do not
need this result.

\section{Hausdorff dimension}\label{s.haus}

\begin{theorem}\label{t.haus}
Suppose that
\begin{equation}\label{e.negmom}
\Eb{W^{-s}}<\infty
\qquad
\text{for all $s\in[0,1)$}
\end{equation}
(in addition to the standard
assumptions
$\Eb{W\log_2W}<1$,
$\Eb{W}=1$,  and $\Pb{W>0}=1$.)
Let $K\subset [0,1]$ be some
(deterministic) nonempty set, let $\Da$ denote its Hausdorff
dimension with respect to the Euclidean metric, and let
$\Db$ denote its Hausdorff dimension with respect to
the random metric $\rho$.
Then a.s.\ $\Db$  is the unique solution of the equation
\begin{equation}
\label{e.Da}
2^{\Da}=\frac{2^{\Db}}{\Eb{W^\Db}}
\end{equation}
in $[0,1]$.
\end{theorem}

As our proof shows, the assumption~\eqref{e.negmom}
may be significantly relaxed.
See Theorems \ref{t.upperbound} and \ref{t.lowerbd}.

\medskip

Now consider the case in which $\log W$ is a Gaussian random variable. Since $\Es{W}=1$,
this implies $W=\exp(\sigma\,Y-\sigma^2/2)$, where $Y$ is a standard Gaussian
of zero mean and unit variance and $\sigma\ge 0$. The assumption~\eqref{e.WlogW}
is then equivalent to the requirement $\sigma^2<\log 4$. In this case,
the moments $\Eb{W^s}$ are easily evaluated and~\eqref{e.Da} gives
\begin{equation}
\label{e.gaussian}
\Da - \Db = \frac{\sigma^2}{\log 4}\,\Db \,(1-\Db)\,, 
\end{equation}
in agreement with~\eqref{e.KPZ}.%
\footnote{In~\eqref{e.KPZ}, $2\,(1-\Daa)$ is the dimension and similarly for $\Dbb$.
The factor of $2$ comes from the fact that the ambient space is two dimensional,
and that the dimension is defined in terms of the measure, not the distance function.
The transition from $\Daa$ to $1-\Daa$ is a passage from the dimension to the
co-dimension, and does not change the form of~\eqref{e.KPZ}.}

For comparison, suppose instead that $W=1\pm\sigma$, each with probability $1/2$.
Then~\eqref{e.WlogW} becomes $|\sigma|<1$ and~\eqref{e.Da} transforms to
$$
2^{\Da}=\frac{2^\Db}{\frac 12 (1-\sigma)^\Db+\frac 12 (1+\sigma)^\Db}\,,
$$
or
$$
\Da = 1+\Db-\log_2\bigl((1-\sigma)^\Db+(1+\sigma)^\Db\bigr).
$$

\bigskip

We now proceed to prove Theorem~\ref{t.haus}.
Define
$$
\phi(s):=s-\log_2 \Eb{W^s}\,.
$$
Then~\eqref{e.Da} reads $\Da=\phi(\Db)$.
The following lemma implies the uniqueness of the $\Db$ satisfying~\eqref{e.Da}.

\begin{lemma}\label{l.psi}
The function $\phi$ is continuous, strictly monotone increasing in $[0,1]$
and maps $[0,1]$ onto $[0,1]$.
\end{lemma}

\proof
Set $\psi(s):= \Eb{(W/2)^s}$.
Continuity of $\psi$ follows from the dominated convergence theorem and convexity
of $\psi$ is immediate by the convexity of $(W/2)^s$ in $s$.
Since $\psi$ is convex and
$$
\psi'(1-)=\Eb{(W/2)\log(W/2)}=\frac 12\, \Eb{W\log W}-\frac 12\log 2<0,
$$
it is strictly monotone decreasing in $[0,1]$.
The lemma follows since $\phi=-\log_2 \psi$, $\phi(0)=0$ and $\phi(1)=1$.
\QED

The following simple lemma can serve to motivate Theorem~\ref{t.haus}, and
is also important in its proof.

\begin{lemma}\label{l.interval}
Let $x,y\in[0,1]$, and let $s\in(0,1]$.
Then
$$
\Eb{\rho(x,y)^s} \le
8\,|x-y|^{\phi(s)}.
$$
\end{lemma}
\proof
Let $[a,b]\in\II_n$.
Then by the construction of $\rho$ and the independence of the different
variables $W_I$, $I\in\II$,
$$
\Eb{\rho(a,b)^s}= 2^{-ns}\,\Eb{W^s}^n\,\Eb{\ell^s}
=|a-b|^{\phi(s)}\,\Eb{\ell^s}.
$$
Now, Jensen's inequality gives $\Eb{\ell^s}\le \Eb{\ell}^s\le 1$.
Note that if $|y-x|\in (2^{-n-1},2^{-n}]$, then
the interval joining $x$ and $y$ can be covered by
two consecutive intervals in $\II_n$, say $[a,b]$ and $[b,c]$.
Then
\begin{equation*}
\begin{aligned}
\Eb{\rho(x,y)^s}
&
\le \Eb{\bigl(\rho(a,b)+\rho(b,c)\bigr)^s}
\\&
\le \Eb{\bigl(2\,\rho(a,b)\bigr)^s+\bigl(2\,\rho(b,c)\bigr)^s}
\\&
= 2^{1+s}\,\Eb{\rho(a,b)^s}
\\&
\le 2^{1+s}\,|a-b|^{\phi(s)}
\le 2^{1+s+\phi(s)}\,|x-y|^{\phi(s)}\,.
\end{aligned}
\end{equation*}
The lemma follows, since $\phi(s)\le 1$ by Lemma~\ref{l.psi}.
\QED

\begin{theorem}
\label{t.upperbound}
Let $K$, $\Da$ and $\Db$ be as in Theorem~\ref{t.haus}.
Then a.s.  $\phi(\Db)\le \Da$.
\end{theorem}

It is worth pointing out that we are not assuming~\eqref{e.negmom} here.

\proof
Let $s\in[0,1]$ and assume that $t:=\phi(s)>\Da$.
We now show that $s \ge \Db$ a.s.
Let $\eps>0$.  Then there is a covering of
$K$ by at most countably many intervals $[x_i,y_i]$
such that $\sum_i |x_i-y_i|^t<\eps$.
By Lemma~\ref{l.interval}, we have
$$
\EB{\sum_i \rho(x_i,y_i)^s} \le
8\, \sum_i |x_i-y_i|^t \le 8\,\eps\,.
$$
By Markov's inequality, with probability at least $1-\sqrt\eps$
we have a covering of $K$ with balls whose radii in the $\rho$ metric
satisfy $\sum r_i^s \le 8\,\sqrt\eps$.
Thus $s\ge \Db$ a.s.
Hence $\Db\le\inf\phi^{-1}(\Da,1]$.
By Lemma~\ref{l.psi}, the theorem follows.
\QED

\begin{theorem}
\label{t.lowerbd}
Let $K$, $\Da$ and $\Db$ be as in Theorem~\ref{t.haus}.
Then a.s.
$$
\Db \ge \sup\bigl\{s\in(0,1): \phi(s)<\Da,\,
\Es{W^{-s}}<\infty\bigr\}.
$$
\end{theorem}

\proof
Suppose that $s\in(0,1)$ satisfies $t:=\phi(s)<\Da$ and $\Es{W^{-s}}<\infty$.
We need to prove that $\Db\ge s$.
Since $\Es{W^s}$ is convex in $s$ and equals to $1$ at $s=0,1$,
we have $t\ge s\ge0$.
Since $\Da>t$, by Frostman's lemma~\cite[Chapter 8]{\Mattila} there is a Borel
probability measure $\nu_0$ supported on $K$
such that
\begin{equation}
\label{e.nuen}
\mathcal E_t(\nu_0):=\iint \frac{d\nu_0(x)\,d\nu_0(y)}{|x-y|^t}<\infty\,.
\end{equation}
Set $a:=\Es{W^s}$, $Z:= W^s/a$ and $Z_I:= W_I^s/a$.
Define
$$
f_n(x):= \prod_{j<n}Z_{I_j(x)},\qquad \nu_n:= f_n\,\nu_0\,.
$$
Since for every $a\in[0,1]$ the sequence $\nu_n[0,a]$ is
a non-negative martingale, it easily follows that
the weak limit $\nu:=\lim_{n\to\infty} \nu_n$ exists.
(See, e.g., the proof of Theorem~\ref{t.bg} in the appendix.)
Then the support of $\nu$ is contained in the
support of $\nu_0$ and therefore in $K$.

Define
$$
\rho_n(x,y):=\rho(x,y)\vee \mu\bigl(I_n(x)\bigr)\vee\mu\bigl(I_n(y)\bigr).
$$
In order to estimate the expectation of
\begin{equation}
\label{e.nunun}
\mathcal E_s(\nu_n;\rho_n):=\iint \frac{d\nu_n(x)\,d\nu_n(y)}{\rho_n(x,y)^s}
\end{equation}
we fix $x,y\in [0,1]$, and estimate
$$
\Eb{f_n(x)\,f_n(y)\,\rho_n(x,y)^{-s}}.
$$
Let $k$ be the smallest integer such that $[x,y]$ contains
some interval of $\II_k$. Then
\begin{equation}
\label{e.xyd}
  |x-y|< 4\cdot 2^{-k}.
\end{equation}
Let $J\in \II_k$ satisfy $J\subset [x,y]$, and let
$J'\in \II_{k-1}$ satisfy $J'\supset J$. Then $x\in J'$ or $y\in J'$.
By symmetry, assume $x\in J'$. Let $\ev G$ denote the $\sigma$-field
$\langle W_{I_j(x)}, W_{I_j(y)}:j< n\rangle$.
Assume first that $k\le n$.
Then
$$
\Eb{\rho(x,y)^{-s}\md \ev G} \le
\Eb{\mu(J)^{-s}\md \ev G} = 2^{ks}\,\Eb{\ell^{-s}}\prod_{j=0}^{k-1} W_{I_j(x)}^{-s}\,.
$$
By Lemma~\ref{l.lowertail} in the appendix and our assumption that 
$\Eb{W^{-s}}<\infty$, we have $\Eb{\ell^{-s}}<\infty$. Using~\eqref{e.xyd},
we therefore obtain
$$
\Eb{\rho(x,y)^{-s}\md \ev G} \le O(1)\,
|x-y|^{-s}\,\prod_{j=0}^{k-1} W_{I_j(x)}^{-s}\,,
$$
where the implied constant may depend on $s$ and the law of $W$.
Now,
$$
f_n(x)=a^{-n}\prod_{j=0}^{n-1} W_{I_j(x)}^s,
$$
and we have a similar expression for $f_n(y)$.
Thus,
$$
\Eb{f_n(x)\,f_n(y)\,\rho(x,y)^{-s}\md \ev G} \le
O(1)\,a^{-2n}\,|x-y|^{-s}\prod_{j=k}^{n-1} W_{I_j(x)}^s\prod_{j=0}^{n-1} W_{I_j(y)}^s\,.
$$
Note that for $j\ge k$ we have $I_j(x)\ne I_j(y)$.
Taking expectations and using the definition of $a$ yields
\begin{multline*}
\Eb{f_n(x)\,f_n(y)\,\rho_n(x,y)^{-s}}\le
\Eb{f_n(x)\,f_n(y)\,\rho(x,y)^{-s}}
\\
\le O(1)\,|x-y|^{-s}\,a^{-k}
\overset{\eqref{e.xyd}}{\le} O(1)\,|x-y|^{-s+\log_2 a}
=
O(1)\,|x-y|^{-t}.
\end{multline*}
Now, if $k>n$, we have instead
$$
\Eb{\rho_n(x,y)^{-s}\md \ev G} \le O(1)\,\Eb{\mu\bigl(I_n(x)\bigr)^{-s}\md\ev G} \le
O(1)\,2^{sn}\prod_{j=0}^{n-1} W_{I_j(x)}^{-s}\,.
$$
The above argument therefore gives in this case,
$$
\Eb{f_n(x)\,f_n(y)\,\rho_n(x,y)^{-s}}
\le O(1)\, 2^{ns}\,a^{-n}\le O(1)\,|x-y|^{-t}.
$$
Thus, we have
$
\Eb{f_n(x)\,f_n(y)\,\rho_n(x,y)^{-s}} \le O(1)\,|x-y|^{-t}
$
for every $x,y\in[0,1]$. Integrating this with respect to $d\nu_0(x)\times d\nu_0(y)$
and applying Fubini, one obtains
\begin{equation}
\label{e.nun}
\Eb{\mathcal E_s(\nu_n;\rho_n)} \le O(1)\,\mathcal E_t(\nu_0)\,.
\end{equation}
Since $\rho_n(x,y)\le \ell$ holds for $x,y\in[0,1]$, this estimate gives
$$
\Eb{(\nu_n[0,1])^2\,\ell^{-s}} \le O(1)\,\mathcal E_t(\nu_0)\,.
$$
Now H\"older's inequality comes into play:
$$
\Eb{(\nu_n[0,1])^{2/(1+s)}} \le
\Eb{(\nu_n[0,1])^2\,\ell^{-s}}^{1/(1+s)} \,\Eb{\ell}^{s/(1+s)}
\le O(1) \,\ev E_t(\nu_0)^{1/(1+s)}.
$$
Thus, the martingale sequence $\nu_n[0,1]$ is uniformly bounded in
$L^p$ with $p=2/(1+s)>1$. It follows by the corresponding martingale
convergence theorem that $\Eb{\nu[0,1]}=\nu_0[0,1]=1$,
and in particular, $\nu[0,1]>0$ with positive probability.
The event $\nu[0,1]>0$ is clearly independent of $\sigma$-field
generated by any finite number of the random variables
$W_{I}$, and therefore has probability $0$ or $1$, and in this case,
$\Pb{\nu[0,1]>0}=1$.

Since a.s.\ $\rho$ is continuous, $\rho_n\to\rho$ uniformly as $n\to\infty$ and $\nu_n\to\nu$
weakly, we have a.s.
$$
\ev E_s(\nu;\rho) \le \liminf_{n\to\infty} \ev E_s(\nu_n;\rho_n)\overset{\eqref{e.nun}}<\infty\,.
$$
The proof is now completed by appealing to
Frostman's criterion~\cite[Chapter 8]{\Mattila},
since $\nu[0,1]>0$ a.s.
\QED

\proofof{Theorem~\ref{t.haus}}
The theorem follows immediately from Lemma~\ref{l.psi} and Theorems~\ref{t.upperbound}
and~\ref{t.lowerbd}.
\QED

\appendix
\section{Some multiplicative cascades background}

\begin{lemma}
\label{l.xlogx}
Our standing assumption $\Eb{W\log_2 W}<1$ implies that $\ell>0$ a.s.
\end{lemma}

\proof
Set $a:=\Eb{W\log_2W}<1$.
We first prove that $\ell>0$ with positive probability.
Since $\ell_n:=\mu_n[0,1]$ is a positive martingale, we have a.s.\
convergence $\ell_n\to \ell$.

The proof will come out of a recurrence relation for
the sequence $b_n:=\Eb{\ell_n\log_2\ell_n}$.
Let $\ell_n'$ and $\ell_n''$ have the law of
$\ell_n$ and be independent and independent
from $W$. Then the law of $\ell_{n+1}$ is the same
as the law of $W\,(\ell_n'+\ell_n'')/2$.
Thus,
$$
\begin{aligned}
&b_{n+1}
=
\Eb{(W/2)\,(\ell_n'+\ell_n'')\log_2 W}+
{}
\\& \qquad\qquad
{}+
\Eb{(W/2)\,(\ell_n'+\ell_n'')\log_2 (\ell_n'+\ell_n'')}+
\Eb{(W/2)\,(\ell_n'+\ell_n'')\log_2 (1/2)}.
\end{aligned}
$$
We now use independence, $\Eb{\ell_n'}=\Eb{W}=1$, the symmetry
between $\ell_n'$ and $\ell_n''$, and the definition of $a$,
and simplify the above to
$$
b_{n+1} = a+ \Eb{\ell_n'\log_2 (\ell_n'+\ell_n'')}-1\,.
$$
Thus,
$$
\begin{aligned}
b_{n+1}-b_n
&
=
a-1+
\Eb{\ell_n'\log_2 (\ell_n'+\ell_n'')}-
\Eb{\ell_n'\log_2 \ell_n'}
\\&
=
a-1+\Eb{\ell_n'\log_2\bigl((\ell_n'+\ell_n'')/\ell_n'\bigr)}.
\end{aligned}
$$
Since $\log_2(1+x/\ell_n')$ is concave as a function of $x$,
and since $\ell_n''$ is independent from $\ell_n'$, Jensen's inequality
applied to the above gives
\begin{multline}
\label{e.bb}
b_{n+1}-b_n\le a-1+\Eb{\ell_n'\log_2(1+\Es{\ell_n''}/\ell_n')}
\\
 =a-1+\Eb{\ell_n\log_2(1+1/\ell_n)}.
\end{multline}
Because $\inf\{b_n:n\in\N\} \ge \inf \{x\log_2 x:x>0\}>-\infty$ and
$a<1$, the set $S:=\{n\in\N: b_{n+1}-b_n>(a-1)/2\}$
is infinite.
For $n\in S$,~\eqref{e.bb} implies
$$
\Eb{g(\ell_n)}>(1-a)/2\,,\qquad \text{where}\qquad
g(x):=x\log_2(1+1/x)\,.
$$
Note that $c:= \sup \{g(x):x>0\}<\infty$.
For all $\eps>0$ and $n\in S$, we have
$$
(1-a)/2 <\Eb{g(\ell_n)} \le
c\,\Pb{\ell_n\ge \eps} + \sup\bigl\{g(x):0<x<\eps\bigr\}\,.
$$
By taking $\eps$ sufficiently small, we can make sure that the last summand
is at most $(1-a)/4$. Then, $ \Pb{\ell_n>\eps}> (1-a)/(4\,c) $.
This proves that $\ell_n$ does not tend to zero in probability, and
hence does not tend to zero a.s.
Thus $\Pb{\ell>0}>0$.

Set $\ell':=2\,\mu[0,1/2]/W_{[0,1]}$ and $\ell'':=2\,\mu[1/2,1]/W_{[0,1]}$.
Then $\ell',\ell''$ and $W_{[0,1]}$ are independent, and each of $\ell'$ and
$\ell''$ has the law of $\ell=(W_{[0,1]}/2)\,(\ell'+\ell'')$.
 But $\ell=0$ if and only if $\ell'=0=\ell''$, since $W>0$ a.s.
Thus,
$$
\Pb{\ell=0}=\Pb{\ell'=0,\,\ell''=0}=\Pb{\ell=0}^2.
$$
 Since $\Pb{\ell>0}>0$, we conclude that $\ell>0$ a.s., which completes the proof.
\QED

\proofof{Theorem~\ref{t.bg}}
Since $\Eb{\mu_n[0,1]}=1$, the sequence of measures $\mu_n$ is tight
in the space of Borel measures on $[0,1]$ with the topology of
weak convergence. Thus, some subsequence
converges to a limit $\mu$.
For every rational $r\in[0,1]\cap\Q$ the sequence $\mu_n[0,r]$ is a positive martingale,
and hence the limit $f(r):=\lim_{n\to\infty} \mu_n[0,r]$ exists
for all $r\in\Q\cap[0,1]$ a.s.
It is immediate to verify that $\mu[0,x)=\sup\bigl\{f(r):r\in\Q\cap[0,x)\bigr\}$
and $\mu[0,1]=f(1)$. This implies that the subsequential limit $\mu$  is
unique, and hence is the weak limit of the entire sequence $\mu_n$.
Thus, the theorem follows from Lemma~\ref{l.xlogx} and the next lemma.
\QED

\begin{lemma}\label{l.noatom}
A.s.\ $\mu$ has no atoms.
\end{lemma}
\proof
Let $Z$ have the distribution of the $\mu$-measure of the
largest atom in $[0,1]$.  Clearly, $\Eb{Z}\le 1$.
Let $Z_1$ and $Z_2$ have the law of
$Z$ with $W,Z_1$ and $Z_2$ independent. Then $W\max(Z_1,Z_2)/2$
has the law of $Z$.  Therefore,
$$
\Eb{Z_1+Z_2}/2 = \Eb{Z}= \Eb{W\max(Z_1,Z_2)/2}=
\Eb{\max(Z_1,Z_2)}/2\,.
$$
Since $Z_1+Z_2\ge \max(Z_1,Z_2)$, and the expectations are the same,
it follows that they are equal a.s.
Thus, on the event $Z_1>0$, we have $Z_2=0$ a.s.
Since $Z_1$ and $Z_2$ are independent, we get $Z_1=0$ a.s.,
as required.
\QED

\begin{lemma}\label{l.lowertail}
Suppose that
\begin{equation}
\label{e.Wr}
\Eb{W^{-r}}<\infty
\end{equation}
for some constant $r>0$.
Then also
$$
\Eb{\ell^{-r}}<\infty\,.
$$
\end{lemma}

A very slightly weaker form of this lemma can be found
in~\cite{MR1847093}, where the assumption~\eqref{e.Wr} is the same,
but the conclusion is that $\Eb{\ell^{-s}}<\infty$ for all $s\in[0,r)$.
However, the setup in~\cite{MR1847093} is more general.

\proof
Let $\ell_1:= 2\,\mu[0,1/2]/W$ and $\ell_2:= 2\,\mu[1/2,1]/W$.
By construction $\ell_1,\ell_2$ and $W$ are independent and each of
$\ell_1$ and $\ell_2$ has
the law of $\ell$.
Moreover
\begin{equation}
\label{e.recursion}
\ell=W\,(\ell_1+\ell_2)/2\,.
\end{equation}
Assume, for a moment, that there is some $\delta>0$
such that
\begin{equation}
\label{e.del}
  \Eb{\ell^{-\delta}}<\infty\,.
\end{equation}
By the means inequality and~\eqref{e.recursion}, we have
\begin{equation}
\label{e.recm}
\ell \ge W \sqrt{\ell_1\,\ell_2}\,.
\end{equation}
Now independence gives for every $s>0$
\begin{equation}
\label{e.boot}
\Eb{\ell^{-s}}\le \Eb{W^{-s}}\,\Eb{\ell_1^{-s/2}}\,\Eb{\ell_1^{-s/2}}
=
\Eb{W^{-s}}\Eb{\ell^{-s/2}}^2.
\end{equation}
Let $S$ denote the set of $s\in[0,\infty)$ such that $\Eb{\ell^{-s}}<\infty$.
By~\eqref{e.del}, we have $\delta\in S$.
Since $\Eb{W^{-s}}$ is a convex function of $s$, the set
$S$ must be an interval. Moreover, $[0,\delta]\subset S$.
Similarly, $\Eb{W^{-s}}<\infty$ for every $s\in[0,r]$.
Now,~\eqref{e.boot} shows that $[0,r]\cap (2 \,S)\subset S$, where
$2\,S=\{2\,s:s\in S\}$.  This implies that $[0,r]\subset S$, as needed.

It remains to prove~\eqref{e.del}. From~\eqref{e.Wr} we know that
\begin{equation}
\label{e.Wtail}
\lim_{x\searrow 0}\, x^{-r}\,\Pb{W<x}=0\,.
\end{equation}
By~\eqref{e.recursion}, for every $b,x>0$, we have
that if $\ell<b\,x/2$ then $W<x$ or $\ell_1+\ell_2<b$.
Thus,
\begin{equation}\label{e.rec1}
\begin{aligned}
\Pb{\ell<b\,x/2}
&
\le
\Pb{W<x}+\Pb{\ell_1+\ell_2<b}
\\&
\le
\Pb{W<x}+\Pb{\ell_1<b,\,\ell_2<b}
\\&
=
\Pb{W<x}+\Pb{\ell<b}^2.
\end{aligned}
\end{equation}
Set $\eps_j:=2^{-2^j-1}$ for $j\in\N$.
Then $\eps_{j+1}=2\,\eps_j^2$.
Let $b_0>0$ satisfy
$$
\Pb{\ell<b_0}\le \eps_0
$$
(Theorem~\ref{t.bg} implies the existence of such a $b_0$).
Set
$$
x_j:=\sup\Bigl\{x:\Pb{W<x}\le \eps_j^2\Bigr\}\,.
$$
Then our assumption~\eqref{e.Wtail} implies
that $x_j\ge \eps_j^{2/r}$ for all
but finitely many $j$.
Inductively define
$$
b_{j+1}=b_j\,x_j/2\,.
$$
Then using induction, the relation $\eps_{j+1}=2\,\eps_j^2$ and
the estimate~\eqref{e.rec1} give
\begin{equation}
\label{e.jj}
\Pb{\ell<b_j}\le \eps_j\,.
\end{equation}
The definition of $b_j$ gives
$$
b_j=2^{-j}\,b_0\,\prod_{k=0}^{j-1} x_k \ge C\, 2^{-j}\prod_{k=0}^{j-1}\eps_k^{2/r}\,,
$$
where $C>0$ is the product of $b_0$ and the finitely many $x_k$
that satisfy $x_k< \eps_k^{2/r}$.
Taking into account the definition of $\eps_j$, we get
$$
b_j\ge C\, 2^{-j-2(j-1+2^j)/r}.
$$
Now,
$$
\Eb{\ell^{-\delta}}
\le
b_0^{-\delta} +
\sum_{j=0}^\infty b_{j+1}^{-\delta}\,\Pb{b_{j+1}\le \ell< b_j}.
$$
For all but finitely many $j$, the above estimate on $b_j$ gives
$b_{j+1}^{-\delta}\le 2^{5 \delta 2^j/r}$, while
$\Pb{b_{j+1}\le \ell< b_j} \le \Pb{\ell< b_j} \overset{\eqref{e.jj}}{\le}
\eps_j =2^{-2^j-1}$.
Thus, we have~\eqref{e.del} for every $\delta\in (0,r/5)$,
and the proof is thus complete.
\QED

\bibliographystyle{halpha}
\bibliography{mr,prep,notmr}

\def\polhk#1{\setbox0=\hbox{#1}{\ooalign{\hidewidth
  \lower1.5ex\hbox{`}\hidewidth\crcr\unhbox0}}} \def\cprime{$'$}
  \def\cprime{$'$}
\begin{thebibliography}{KPZ88}

\bibitem[DS08]{DSinprep}
Bertrand Duplantier and Scott Sheffield, 2008.
\newblock In preparation.

\bibitem[Kol91]{MR1124922}
A.~N. Kolmogorov.
\newblock The local structure of turbulence in incompressible viscous fluid for
  very large {R}eynolds numbers.
\newblock {\em Proc. Roy. Soc. London Ser. A}, 434(1890):9--13, 1991.
\newblock Translated from the Russian by V. Levin, Turbulence and stochastic
  processes: Kolmogorov's ideas 50 years on.

\bibitem[KP76]{MR0431355}
J.-P. Kahane and J.~Peyri{\`e}re.
\newblock Sur certaines martingales de {B}enoit {M}andelbrot.
\newblock {\em Advances in Math.}, 22(2):131--145, 1976.

\bibitem[KPZ88]{MR947880}
V.~G. Knizhnik, A.~M. Polyakov, and A.~B. Zamolodchikov.
\newblock Fractal structure of {$2$}{D}-quantum gravity.
\newblock {\em Modern Phys. Lett. A}, 3(8):819--826, 1988.

\bibitem[Liu01]{MR1847093}
Quansheng Liu.
\newblock Asymptotic properties and absolute continuity of laws stable by
  random weighted mean.
\newblock {\em Stochastic Process. Appl.}, 95(1):83--107, 2001.

\bibitem[Man74]{MandelbrotCascades}
B.~Mandelbrot.
\newblock Intermittent turbulence in self similar cascades: divergence of high
  moments and dimension of carrier.
\newblock {\em J. Fluid Mech.}, 62:331--333, 1974.

\bibitem[Mat95]{MR1333890}
Pertti Mattila.
\newblock {\em Geometry of sets and measures in {E}uclidean spaces}, volume~44
  of {\em Cambridge Studies in Advanced Mathematics}.
\newblock Cambridge University Press, Cambridge, 1995.
\newblock Fractals and rectifiability.

\bibitem[OW00]{MR1835030}
Mina Ossiander and Edward~C. Waymire.
\newblock Statistical estimation for multiplicative cascades.
\newblock {\em Ann. Statist.}, 28(6):1533--1560, 2000.

\bibitem[Pol87]{MR913671}
A.~M. Polyakov.
\newblock Quantum gravity in two dimensions.
\newblock {\em Modern Phys. Lett. A}, 2(11):893--898, 1987.

\bibitem[Yag66]{YaglomTurbulence}
A.M. Yaglom.
\newblock Effect of fluctuations in energy dissipation rate on the form of
  turbulence characteristics in the inertial subrange.
\newblock {\em Dokl. Akad. Nauk SSSR}, 166:49--52, 1966.

\end{thebibliography}

\end{document}